\documentclass[11pt]{article}
\usepackage[german,english]{babel}
\usepackage{amsmath}
\usepackage{amssymb}
\usepackage{amsthm}

\usepackage{url}

\newcommand\dstyle\displaystyle

\newcommand\mLP{\\[\medskipamount]}
\newcommand\mPP{\\[\medskipamount]\indent}

\newcommand\RR{\mathbb{R}}
\newcommand\ZZ{\mathbb{Z}}
\newcommand\FSJ{{\cal J}}
\newcommand\al\alpha
\newcommand\be\beta
\newcommand\ga\gamma
\newcommand\de\delta
\newcommand\tha\theta
\newcommand\la\lambda
\newcommand\si\sigma
\newcommand\Ga{\Gamma}

\newcommand\half{\frac12}
\newcommand\thalf{\tfrac12}
\newcommand\pa\partial
\newcommand\iy\infty
\newcommand\wt{\widetilde}
\newcommand\wh{\widehat}
\newcommand\const{\mathrm{const.}\,}

\newcommand{\hyp}[5]{\,\mbox{}_{#1}F_{#2}\!\left(
  \genfrac{}{}{0pt}{}{#3}{#4};#5\right)}
\newcommand{\qhyp}[5]{\,\mbox{}_{#1}\phi_{#2}\!\left(
  \genfrac{}{}{0pt}{}{#3}{#4};#5\right)}
\newcommand{\bisub}[2]{\genfrac{}{}{0pt}{1}{#1}{#2}}

\numberwithin{equation}{section}

\begin{document}
\title{A short survey of duality in special functions}
\author{Tom H. Koornwinder}
\date{}
\maketitle
\begin{abstract}
This is a tutorial on duality properties of special functions, mainly of
orthogonal polynomials in the ($q$-)Askey scheme.
It is based on the first part of the 2017
R.~P. Agarwal Memorial Lecture delivered by the author.
\end{abstract}
\section{Introduction}
Classical orthogonal polynomials $p_n(x)$ and their generalizations in
the Askey and $q$-Askey scheme have the property that they are eigenfunctions
of some second order operator $L$ with eigenvalues depending on $n$, which
therefore may be called the spectral variable. Moreover, being orthogonal
polynomials, the $p_n(x)$ satisfy a three-term recurrence relation and are
therefore, as functions of $n$, eigenfunctions of a so-called Jacobi operator
with eigenvalues $x$. This duality phenomenon was also guiding for the
author in the companion paper \cite{33}, where he derived the dual
addition formula for continuous $q$-ultraspherical polynomials.

This paper gives a brief tutorial type survey of duality, mainly for orthogonal
polynomials, but also a little bit for transcendental special functions.
This paper is based on the first part of the
\emph{R.~P. Agarwal Memorial Lecture}, which
the author delivered on November 2, 2017 during the conference
ICSFA-2017 held in Bikaner, Rajasthan, India. See \cite{33} for the paper
based on the second part.

With pleasure I
remember to have met Prof.\ Agarwal during
the workshop on Special Functions and Differential Equations
held at the Institute of Mathematical Sciences in Chennai,
January 1997, where he delivered the opening address \cite{14}.
I cannot resist to quote from it the following wise words, close
to the end of the article:
\mLP
``{\sl 
I think that I have taken enough time and I close my discourse- with a word of caution and advice to the research workers in the area of special functions and also those who use them in physical problems. The corner stones of classical analysis are
`elegance, simplicity, beauty and perfection.' Let them not be lost in your work. Any analytical generalization of a special function, only for the sake of a generalization by adding a few terms or parameters here and there, leads us nowhere. All research work should be meaningful and aim at developing a quality technique or have a bearing in some allied discipline.}''
\paragraph{Note}
For definition and notation of ($q$-)shifted factorials and
($q$-)hypergeometric series see \cite[\S1.2]{3}.
In the $q=1$ case we will mostly meet terminating hypergeometric series
\begin{equation}
\hyp rs{-n,a_2,\ldots,a_r}{b_1,\ldots,b_s}z:=
\sum_{k=0}^n \frac{(-n)_k}{k!}\,\frac{(a_2,\ldots,a_r)_k}
{(b_1,\ldots,b_s)_k}\,z^k.
\label{76}
\end{equation}
Here $(b_1,\ldots,b_s)_k:=(b_1)_k\ldots(b_s)_k$ and
$(b)_k:=b(b+1)\ldots(b+k-1)$ is the \emph{Pochhammer symbol} or
\emph{shifted factorial}. In \eqref{76} we even allow that
$b_i=-N$ for some $i$ with $N$ integer $\ge n$. There is no problem
because the sum on the right terminates at $k=n\le N$.

In the $q$-case we will always assume that $0<q<1$.
We will only meet terminating $q$-hypergeometric series of the form
\begin{equation}
\qhyp{s+1}s{q^{-n},a_2,\ldots,a_{s+1}}{b_1,\ldots,b_s}{q,z}:=
\sum_{k=0}^n \frac{(q^{-n};q)_k}{(q;q)_k}\,
\frac{(a_2,\ldots,a_{s+1};q)_k}{(b_1,\ldots,b_s;q)_k}\,z^k.
\label{77}
\end{equation}
Here $(b_1,\ldots,b_s;q)_k:=(b_1;q)_k\ldots(b_s;q)_k$ and
$(b;q)_k:=(1-b)(1-qb)\ldots(1-q^{k-1}b)$ is the
\emph{$q$-Pochhammer symbol} or \emph{$q$-shifted factorial}.
In \eqref{77} we even allow that $b_i=q^{-N}$ for some $i$ with
$N$ integer $\ge n$.

For formulas on orthogonal polynomials in the ($q$-)Askey scheme we
will often refer to Chapters 9 and 14 in \cite{2}. Almost all
of these formulas, with different numbering, are available in open
access on \url{http://aw.twi.tudelft.nl/~koekoek/askey/} .
\section{The notion of duality in special functions}
\label{72}
With respect to a (positive) measure $\mu$ on $\RR$ with support
containing infinitely many
points we can define \emph{orthogonal polynomials} (OPs)
$p_n$ ($n=0,1,2,\ldots$), unique up to nonzero real constant factors,
as (real-valued) polynomials $p_n$ of degree $n$ such that
\begin{equation*}
\int_\RR p_m(x)\,p_n(x)\,d\mu(x)=0\qquad(m,n\ne0).
\end{equation*}
Then the polynomials $p_n$ satisfy a \emph{three-term recurrence relation}
\begin{equation}
x\,p_n(x)=A_n\,p_{n+1}(x)+B_n\,p_n(x)+C_n\,p_{n-1}(x)
\qquad(n=0,1,2,\ldots),
\label{52}
\end{equation}
where the term $C_n\,p_{n-1}(x)$ is omitted if $n=0$, and where
$A_n,B_n,C_n$ are real and
\begin{equation}
A_{n-1}C_n>0\qquad(n=1,2,\ldots).
\label{53}
\end{equation}
By \emph{Favard's theorem} \cite{15}
we can conversely say that if $p_0(x)$ is
a nonzero real constant, and the
$p_n(x)$ ($n=0,1,2,\ldots$) are generated
by \eqref{52}
for certain real $A_n,B_n,C_n$ which satisfy \eqref{53}, then the $p_n$
are OPs with respect to a certain measure $\mu$ on
$\RR$.

With $A_n,B_n,C_n$ as in \eqref{52} define a \emph{Jacobi operator} $M$,
acting on infinite sequences $\{g(n)\}_{n=0}^\iy$, by
\begin{equation*}
(Mg)(n)=M_n\big(g(n)\big):=A_n\,g(n+1)+B_n\,g(n)+C_n\,g(n-1)
\qquad(n=0,1,2,\ldots),
\end{equation*}
where the term $C_n\,g(n-1)$ is omitted if $n=0$. Then \eqref{52}
can be rewritten as the eigenvalue equation
\begin{equation}
M_n\big(p_n(x)\big)=x\,p_n(x)\qquad(n=0,1,2,\ldots).
\label{54}
\end{equation}
One might say that the study of a system of OPs $p_n$
turns down to the spectral theory and harmonic analysis associated
with the operator $M$. From this perspective one can wonder if
the polynomials $p_n$ satisfy some \emph{dual} eigenvalue equation
\begin{equation}
(Lp_n)(x)=\la_n\,p_n(x)
\label{55}
\end{equation}
for $n=0,1,2,\ldots$,
where $L$ is some linear operator acting on the space of polynomials.
We will consider varioua types of operators $L$ together with the
corresponding OPs, first in the Askey scheme and next in the
$q$-Askey scheme.
\subsection{The Askey scheme}
\label{71}
\paragraph{Classical OPs}
\emph{Bochner's theorem} \cite{16}
classifies the second order differentai operators
$L$ together with the OPs $p_n$ such that
\eqref{55} holds for certain eigenvalues $\la_n$.
The resulting \emph{classical orthogonal polynomials} are essentially
the polynomials listed in the table below. Here $d\mu(x)=w(x)\,dx$
on $(a,b)$ and the closure of that interval is the support of $\mu$.
Furthermore, $w_1(x)$ occurs in the formula for $L$ to be given after
the table.
\begin{center}
\begin{tabular}{|c|c|c|c|c|c|c|}
\hline
name&$p_n(x)$&$w(x)$&$\frac{w_1(x)}{w(x)}$&$(a,b)$&constraints&$\la_n$\\
\hline
Jacobi&$P_n^{(\al,\be)}(x)$&{\small$(1-x)^\al(1+x)^\be$}&{\small$1-x^2$}&$(-1,1)$&$\al,\be>-1$&{\small$-n(n+\al+\be+1)$}\\
Laguerre&$L_n^{(\al)}(x)$&$x^\al e^{-x}$&$x$&$(0,\iy)$&$\al>-1$&$-n$\\
Hermite&$H_n(x)$&$e^{-x^2}$&$1$&$(-\iy,\iy)$&&$-2n$\\
\hline
\end{tabular}
\end{center}
Then
\begin{equation*}
(Lf)(x)=w(x)^{-1} \frac d{dx}\big( w_1(x)\,f'(x)\big).
\end{equation*}

For these classical OPs the duality goes much further than the two
dual eigenvalue equations \eqref{54} and \eqref{55}.
In particular for Jacobi polynomials it is true to a large extent
that every formula or property involving $n$ and $x$ has a dual
formula or property where the roles of $n$ and $x$ are interchanged.
We call this the \emph{duality principle}.
If the partner formula or property is not yet known then it is usually
a good open problem to find it (but one should be warned that there
are examples where the duality fails).

The Jacobi, Laguerre and Hermite families
are connected by limit transitions, as is already suggested by limit
transitions for their (rescaled) weight functions:
\begin{itemize}
\item
Jacobi $\to$ Laguerre:\quad $x^\al(1-\be^{-1}x)^\be\to x^\al e^{-x}$\quad
as\quad$\be\to\iy$;
\item
Jacobi $\to$ Hermite:\quad $\big(1-\al^{-1}x^2\big)^\al\to e^{-x^2}$\quad
as\quad$\al\to\iy$;
\item
Laguerre $\to$ Hermite:\quad
$e^{\al(1-\log\al)}\big((2\al)^\half x+\al\big)^\al
e^{-(2\al)^\half x-\al}\to e^{-x^2}$\quad as\quad$\al\to\iy$.
\end{itemize}
Formulas and properties of the three families can be expected to be
connected under these limits. Although this is not always the case,
this \emph{limit principle} is again a good source of open problems.
\paragraph{Discrete analogues of classical OPs}
Let $L$ be a second order difference operator:
\begin{equation}
(Lf)(x):=a(x)\,f(x+1)+b(x)\,f(x)+c(x)\,f(x-1).
\label{56}
\end{equation}
Here as solutions of \eqref{55} we will also allow OPs
$\{p_n\}_{n=0}^N$ for some \emph{finite} $N\ge0$, which will be orthogonal
with respect to positive weights $w_k$ ($k=0,1,\ldots,N$) on a finite
set of points $x_k$ ($k=0,1,\ldots,N$):
\begin{equation*}
\sum_{k=0}^N p_m(x_k)\,p_n(x_k)\,w_k=0\qquad(m,n=0,1,\ldots,N;\;m\ne n).
\end{equation*}
If such a finite system of OPs satisfies \eqref{55} for $n=0,1,\ldots,N$
with $L$ of the form \eqref{56} then the highest $n$ for which
the recurrence relation \eqref{52}
holds is $n=N$, where the zeros of $p_{N+1}$ are precisely the
points $x_0,x_1,\ldots,x_N$.

The classification of OPs satisfying \eqref{55}
with $L$ of the form \eqref{56} (first done by O. Lancaster, 1941,
see \cite{17})
yields the four families of
Hahn, Krawtchouk, Meixner and Charlier polynomials, of which
Hahn and Krawtchouk are finite systems, and Meixner and Charlier
infinite systems with respect to weights on countably infinite sets.

\emph{Krawtchouk polynomials} \cite[(9.11.1)]{2} are given by
\begin{equation}
K_n(x;p,N):=
\hyp21{-n,-x}{-N}{p^{-1}}\quad(n=0,1,2,\ldots,N).
\label{57}
\end{equation}
They satsify the orthogonality relation
\begin{equation*}
\sum_{x=0}^N (K_m\,K_n\,w)(x;p,N)=\frac{(1-p)^N}{w(n;p,N)}\,\de_{m,n}
\end{equation*}
with weights
\begin{equation*}
w(x;p,N):=\binom Nx p^x(1-p)^{N-x}\quad(0<p<1).
\end{equation*}
By \eqref{57} they are \emph{self-dual}\,:
\begin{equation*}
K_n(x;p,N)=K_x(n;p,N)\qquad(n,x=0,1,\ldots,N).
\end{equation*}
The three-term recurrence relation \eqref{54} immediately
implies a dual equation \eqref{55} for such OPs.

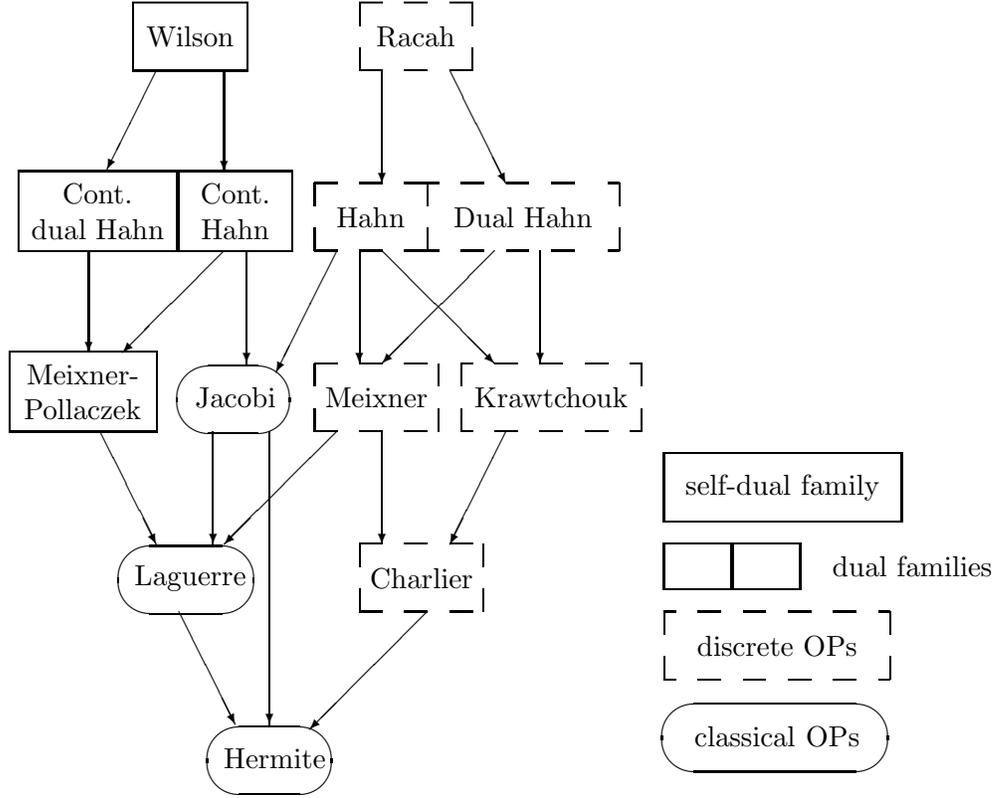
\begin{figure}[t]
\setlength{\unitlength}{3mm}
\begin{picture}(30,36)
\put(8.5,34.5) {\framebox(5,3) {Wilson}}
\put(9.5,34.5) {\vector(-1,-2){2.2}}
\put(12.5,34.5) {\vector(0,-1){4.5}}
\put(18.5,34.5) {\dashbox(5,3) {Racah}}
\put(19.5,34.5) {\vector(0,-1){5}}
\put(22.5,34,5) {\vector(1,-2){2.5}}
\put(3.4,26.5) {\framebox(7,3.5) {\shortstack {Cont.\\[1mm]dual Hahn}}}
\put(6.5,26.5) {\vector(0,-1){4.5}}
\put(10.5,26.5) {\framebox(5,3.5) {\shortstack {Cont.\\[1mm]Hahn}}}
\put(12.5,26.5) {\vector(-1,-1){4.5}}
\put(13.5,26.5) {\vector(0,-1){5}}
\put(16.5,26.5) {\dashbox(5,3) {Hahn}}
\put(17.5,26.5) {\vector(-1,-2){2.7}}
\put(18.5,26.5) {\vector(0,-1){5}}
\put(19.5,26.5) {\vector(1,-1){5}}
\put(21.5,26.5) {\dashbox(8.5,3) {Dual Hahn}}
\put(24.5,26.5) {\vector(-1,-1){5}}
\put(26.5,26.5) {\vector(0,-1){5}}
\put(3,18.5) {\framebox(6.5,3.5)
{\shortstack {Meixner-\\[1mm]Pollaczek}}}
\put(7,18.5) {\vector(1,-2){2.5}}
\put(12.9,19.9) {\oval(5,3)}
\put(10.5,18.5) {\makebox(5,3) {Jacobi}}
\put(12.0,18.5) {\vector(0,-1){5}}
\put(14.5,18.5) {\vector(0,-1){13}}
\put(16.5,18.5) {\dashbox(5.5,3) {Meixner}}
\put(17.5,18.5) {\vector(-1,-1){5}}
\put(19.5,18.5) {\vector(0,-1){5}}
\put(23,18.5) {\dashbox(8,3) {Krawtchouk}}
\put(25,18.5) {\vector(-1,-2){2.5}}
\put(10.8,11.9) {\oval(6,3)}
\put(8,10.5) {\makebox(6,3) {Laguerre}}
\put(10.5,10.5) {\vector(1,-2){2.5}}
\put(18.5,10.5) {\dashbox(5.5,3) {Charlier}}
\put(21.5,10.5) {\vector(-1,-1){5.2}}
\put(14.5,3.9) {\oval(5.5,3)}
\put(12.0,2.5) {\makebox(5.5,3) {Hermite}}
\put(32,14.5) {\framebox(10.5,3) {self-dual family}}
\put(32,11.5) {\framebox(3,2) {}}
\put(35,11.5) {\framebox(3,2) {}}
\put(39,11.5) {\makebox(8,2) {dual families}}
\put(32,7.5) {\dashbox(10,3) {discrete OPs}}
\put(36.9,4.9) {\oval(10,3)}
\put(32,3.5) {\makebox(10,3) {classical OPs}}
\end{picture}
\caption{The Askey scheme}
\label{fig:1}
\end{figure}

The four just mentioned families of discrete OPs
are also connected by limit relations. Moreover, the classical OPs
can be obtained as limit cases of them, but not conversely.
For instance, \emph{Hahn polynomials} \cite[(9.5.1)]{2}
are given by
\begin{equation}
Q_n(x;\al,\be,N):=\hyp32{-n,n+\al+\be+1,-x}{\al+1,-N}1\qquad
(n=0,1,\ldots,N)
\label{58}
\end{equation}
and they satisfy the orthogonality relation
\begin{equation*}
\sum_{x=0}^N (Q_mQ_n w)(x;\al,\be,N)=0\qquad(m,n=0,1,\ldots,N;\;m\ne n;
\;\al,\be>-1)
\end{equation*}
with weights
\begin{equation*}
w(x;\al,\be,N):=\frac{(\al+1)_x\,(\be+1)_{N-x}}{x!\,(N-x)!}\,.
\end{equation*}
Then by \eqref{58} (rescaled) Hahn polynomials tend to (shifted) Jacobi
polynomials:
\begin{equation}
\lim_{N\to\iy}Q_n(Nx;\al,\be,N)=
\hyp21{-n,n+\al+\be+1}{\al+1}x=
\frac{P_n^{(\al,\be)}(1-2x)}{P_n^{(\al,\be)}(1)}\,.
\label{66}
\end{equation}
\paragraph{Continuous versions of Hahn and Meixner polynomials}\quad\\
A variant of the difference operator \eqref{56} is the operator
\begin{equation}
(Lf)(x):=A(x)\,f(x+i)+B(x)\,f(x)+\overline{A(x)}\,f(x-i)\qquad(x\in\RR),
\label{59}
\end{equation}
where $B(x)$ is real-valued.
Then further OPs satisfying \eqref{55} are
the continuous Hahn polynomials and the Meixner-Pollaczek polynomials
\cite[Ch.~9]{2}.
\paragraph{Insertion of a quadratic argument}\quad\\
For an operator $\wt L$ and some polynomial $\si$
of degree 2 we can define an operator $L$ by
\begin{equation}
(Lf)\big(\si(x)\big):=\wt{L}_x\Big(f\big(\si(x)\big)\Big),
\label{61}
\end{equation}
Now we look for OPs satisfying \eqref{55} where $\wt L$ is
of type \eqref{56} or \eqref{59}. So
\begin{equation}
\wt L_x\Big(p_n\big(\si(x)\big)\Big)=\la_n\,p_n\big(\si(x)\big).
\label{60}
\end{equation}
The resulting OPs are the
Racah polynomials and dual Hahn polynomials for \eqref{60} with
$\wt L$ of type \eqref{56}, and Wilson polynomials and continuous
dual Hahn polynomials for \eqref{60} with
$\wt L$ of type~\eqref{59}, see again \cite[Ch.~9]{2}.
\mPP
The OPs satisfying \eqref{55} in the cases discussed until now
form together the \emph{Askey scheme}, see Figure \ref{fig:1}.
The arrows denote limit transitions.

In the Askey scheme we emphasize the self-dual families:
Racah, Meixner, Krawtchouk and Charlier for the OPs with discrete
orthogonality measure,
and Wilson and Meixner-Pollaczek for the OPs with non-discrete
orthogonality measure. We already met perfect self-duality for
the Krawtchouk polynomials, which is also the case for Meixner and
Charlier polynomials. For the Racah polynomials the dual OPs are
still Racah polynomials, but with different values of the parameters:
\begin{multline*}
R_n\big(x(x+\de-N);\al,\be,-N-1,\de\big)
:=
\hyp43{-n,n+\al+\be+1,-x,x+\de-N}{\al+1,\be+\de+1,-N}1\\
=R_x(n(n+\al+\be+1);-N-1,\de,\al,\be)\qquad
(n,x=0,1,\ldots,N).
\end{multline*}
The orthogonality relation for these Racah polynomials involves a weighted
sum of terms $(R_mR_n)\big(x(x+\de-N);\al,\be,-N-1,\de\big)$ over
$x=0,1,\ldots N$.

For Wilson polynomials we have also self-duality with a change of
parameters but the self-duality is not perfect, i.e., not related to the
orthogonality relation:
\begin{multline}
\const {W_n(x^2;a,b,c,d)}:=
\hyp43{-n,n+a+b+c+d-1,a+ix,a-ix}{a+b,a+c,a+d}1\\
=\const W_{-ix-a}\Big(\big(i(n+a')\big)^2;a',b',c',d'\Big),
\label{62}
\end{multline}
where $a'=\thalf(a+b+c+d-1)$, $a'+b'=a+b$, $a'+c'=a+c$, $a'+d'=a+d$.
The duality \eqref{62} holds for $-ix-a=0,1,2,\ldots$, while the orthogonality
relation for the Wilson polynomials involves a weighted integral
of $(W_mW_n)(x^2;a,b,c,d)$ over $x\in[0,\iy)$.

As indicated in Figure \ref{fig:1}, the dual Hahn polynomials
\begin{equation*}
R_n\big(x(x+\al+\be+1);\al,\be,N\big):=
\hyp32{-n,-x,x+\al+\be+1}{\al+1,-N}1\qquad(n=0,1,\ldots,N)
\end{equation*}
are dual to the Hahn polynomials \eqref{58}:
\begin{equation*}
Q_n(x;\al,\be,N)=R_x\big(n(n+\al+\be+1);\al,\be,N\big)\qquad
(n,x=0,1,\ldots,N).
\end{equation*}
The duality is perfect: the dual orthogonality relation for the Hahn
polynomials is the orthogonality relation for the dual Hahn polynomials,
and conversely.
There is a similar, but non-perfect duality between continuous Hahn
and continuous dual Hahn.

The classical OPs are in two senses exceptional within the
Askey scheme. First, they are the only families which are not
self-dual or dual to another family of OPs.
Second, they are the only
continuous families which are not related by analytic continuation
to a discrete family.

With the arrows in the Askey scheme given it can be taken as
a leading principle to link also the formulas and properties of
the various families in the Askey scheme by these arrows.
In particular, if one has some formula or property for a family lower
in the Askey scheme, say for Jacobi, then one may look for  the
corresponding formula or property higher up, and try to find it if it
is not yet known. In particular, if one could find the result on the highest
Racah or Wilson level, which is self-dual then, going down along the
arrows, one might also obtain two mutually dual results in the Jacobi case.
\subsection{The $q$-Askey scheme}
The families of OPs in the $q$-Askey scheme\footnote{See
\url{http://homepage.tudelft.nl/11r49/pictures/large/q-AskeyScheme.jpg}}
\cite[Ch.~14]{2}
result from the classification \cite{21}, \cite{18}, \cite{19},
\cite{20}
of OPs satisfying \eqref{55}, where $L$ is defined in terms of
the operator $\wt L$ and the function $\si$ by \eqref{61}, where
$\wt L$ is of type \eqref{56} or \eqref{59}, and where
$\si(x)=q^x$ or equal to a quadratic polynomial in $q^x$. This choice
of $\si(x)$ is the new feature deviating from what we discussed about the
Askey scheme. And here $q$ enters, with $0<q<1$ always assumed.
The $q$-Askey scheme is considerably larger than the Askey scheme,
but many features of the Askey scheme return here, in particular it has
arrows denoting limit relations. Moreover, the $q$-Askey scheme is
quite parallel to the Askey scheme in the sense that
OPs in th $q$-Askey scheme,
after suitable rescaling, tend to OPs in the Askey scheme as
$q\uparrow1$. Parallel to Wilson and Racah polynomials at the top
of the Askey scheme there are Askey--Wilson polynomials \cite{9}
and $q$-Racah polynomials at the top of the $q$-Askey scheme.
These are again self-dual families, with the self-duality for
$q$-Racah being perfect.

The guiding principles discussed before about formulas or properties
related  by duality or limit transitions now extend to the $q$-Askey scheme:
both within the $q$-Askey scheme and in relation to the Askey scheme
by letting $q\uparrow1$. For instance, one can hope to find as many
dual pairs of significant formulas and properties of Askey--Wilson
polynomials as
possible which have mutually dual limit cases for Jacobi
polynomials. In fact, we realize this in \cite{33}
with the addition and dual addition formula by taking limits
from the continuous $q$-ultraspherical polynomials (a self-dual
one-parameter subclass of the four-parameter class of
Askey--Wilson polynomials) to the ultraspherical polynomials
(a one-parameter subclass of the two-parameter class of Jacobi
polynomials).

One remarkable aspect of duality in the two schemes concerns
the discrete OPs living there.
Leonard (1982)
classified all systems of OPs $p_n(x)$ with respect to weights
on a countable set $\{x(m)\}$
for which there is a system
of OPs $q_m(y)$ on a countable set $\{y(n)\}$
such that
\[
p_n\big(x(m)\big)=q_m\big(y(n)\big).
\]
His classification yields the OPs in the $q$-Askey scheme
which are orthogonal with respect to weights on a countable set
together with their limit cases for $q\uparrow1$ and $q\downarrow-1$
(where we allow $-1<q<1$ in the $q$-Askey scheme).
The $q\downarrow -1$ limit case yields the
Bannai--Ito polynomials~\cite{22}.
\subsection{Duality for non-polynomial special functions}
For Bessel functions $J_\al$ see \cite[Ch.~10]{11}
and references given there.
It is convenient to use a different standardization and notation:
\begin{equation*}
\FSJ_\al(x):=\Ga(\al+1)\,(2/x)^\al\,J_\al(x).
\end{equation*}
Then (see \cite[(10.16.9)]{11})
\[
\FSJ_\al(x)=
\sum_{k=0}^\iy\frac{(-\tfrac14 x^2)^k}{(\al+1)_k\,k!}
=\hyp01-{\al+1}{-\tfrac14x^2}\qquad(\al>-1).
\]
$\FSJ_\al$ is an even entire analytic function. Some special cases are
\begin{equation}
\FSJ_{-1/2}(x)=\cos x,\quad
\FSJ_{1/2}(x)=\frac{\sin x}x\,.
\label{63}
\end{equation}
The \emph{Hankel transform} pair \cite[\S10.22(v)]{11}, for $f$ in a suitable
function class, is given by
\begin{equation*}
\begin{cases}
&\dstyle\wh f(\la)=\int_0^\iy f(x)\FSJ_\al(\la x) x^{2\al+1}\,dx,\mLP
&\dstyle
f(x)=\frac1{2^{2\al+1}\Ga(\al+1)^2}\int_0^\iy \wh f(\la)
\FSJ_\al(\la x) \la^{2\al+1}\,d\la.
\end{cases}
\end{equation*}
In view of \eqref{63} the Hankel transform contains the
Fourier-cosine and Fourier-sine transform as special cases for
$\al=\pm\half$.

The functions $x\mapsto\FSJ_\al(\la x)$ satisfy the eigenvalue equation
\cite[(10.13.5)]{11}
\begin{equation}
\left(\frac{\pa^2}{\pa x^2}+\frac{2\al+1}x\,\frac \pa{\pa x}\right)
\FSJ_\al(\la x)=-\la^2\,\FSJ_\al(\la x).
\label{64}
\end{equation}
Obviously, then also
\begin{equation}
\left(\frac{\pa^2}{\pa \la^2}+\frac{2\al+1}\la\,\frac \pa{\pa\la}\right)
\FSJ_\al(\la x)=-x^2\,\FSJ_\al(\la x).
\label{65}
\end{equation}
The differential operator in \eqref{65} involves the
spectral variable $\la$
of \eqref{64}, while the eigenvalue in \eqref{65} involves the
$x$-variable in the differential operator in \eqref{64}.

The Bessel functions and the Hankel transform are
closely related to the Jacobi polynomials \eqref{66}
and their orthogonality relation. Indeed, we have the limit formulas
\begin{equation*}
\lim_{n\to\iy}\frac{P_n^{(\al,\be)}\big(\cos(n^{-1}x)\big)}
{P_n^{(\al,\be)}(1)}=\FSJ_\al(x),\qquad
\lim_{\bisub{\vphantom{|}\nu\to\iy}{\vphantom{|}\nu\la=1,2,\ldots}}
\frac{P_n^{(\al,\be)}\big(\cos(\nu^{-1}x)\big)}
{P_n^{(\al,\be)}(1)}=\FSJ_\al(\la x).
\end{equation*}
There are many other examples of non-polynomial special functions
being limit cases of OPs in the ($q$-)Askey scheme, see for instance
\cite{23}, \cite{24}.

In 1986 Duistermaat \& Gr\"unbaum \cite{12} posed the question
if the pair of eigenvalue equations \eqref{64}, \eqref{65}
could be generalized to a pair
\begin{equation}
\begin{split}
L_x\big(\phi_\la(x)\big)&=-\la^2\,\phi_\la(x),\\
M_\la\big(\phi_\la(x)\big)&=\tau(x)\,\phi_\la(x)
\end{split}
\label{67}
\end{equation}
for suitable differential operators $L_x$ in $x$ and $M_\la$ in $\la$
and suitable functions $\phi_\la(x)$ solving the two equations.
Here the functions $\phi_\la(x)$ occur as eigenfunctions in two ways:
for the operator $L_x$ with eigenvalue depending on $\la$ and for
the operator $M_\la$ with eigenvalue depending on $x$.
Since the occurring eigenvalues of an operator form its spectrum,
a phenomenon as in \eqref{67} is called \emph{bispectrality}.
For the case of a second order differential operator $L_x$
written in potential form $L_x=d^2/dx^2-V(x)$ they classified
all possibilities for \eqref{67}. Beside the mentioned Bessel cases
and a case with Airy functions (closely related to Bessel functions)
they obtained two other families where $M_\la$ is a higher than
second order differential operator. These could be obtained by
successive \emph{Darboux transformations} applied to $L_x$
in potential form.
A Darboux transformation produces
a new potential from a given potential $V(x)$
by a formula which involves an eigenfunction of $L_x$ with eigenvalue 0.
Their two new families get a start by the application of
a Darboux transformation
to the Bessel differential equation \eqref{64}, rewritten in
potential form
\begin{equation*}
\phi_\la''(x)-(\al^2-\tfrac14)x^{-2}\phi_\la(x)=-\la^2\phi_\la(x),\qquad
\phi_\la(x)=(\la x)^{\al+\half} \FSJ_\al(\la x).
\end{equation*}
Here $\al$ should be in $\ZZ+\thalf$ for a start of the first new family
or in $\ZZ$ for a start of the second new family. For other values
of $\al$ one would not obtain a dual eigenvalue equation with
$M_\la$ a finite order differential operator.

Just as higher order differential operators $M_\la$ occur in \eqref{67},
there has been a lot of work on studying OPs satisfying
\eqref{55} with $L$ a higher order differential operator.
See a classification in \cite{25}, \cite{26}. All occurring
OPs, the so-called \emph{Jacobi type} and
\emph{Laguerre type polynomials},
have a Jacobi or Laguerre orthogonality measure with integer
values of the parameters, supplemented by mass points at one or
both endpoints of the orthogonality interval. Some of the
Bessel type functions in the second new class in \cite{12} were
obtained in \cite{27} as limit cases of Laguerre type polynomials.
\subsection{Some further cases of duality}
The self-duality property of the family of Askey-Wikson polynomials
is reflected in Zhedanov's \emph{Askey--Wilson algebra} \cite{28}.
A larger algebraic structure is the \emph{double affine Hecke
algebra} (DAHA), introduced by Cherednik and extended by Sahi.
The related special functions are so-called
\emph{non-symmetric} special functions. They are functions in several
variables and associated with root systems. Again there is a duality,
both in the DAHA and for the related special functions.
For the (one-variable) case
of the non-symmetric Askey--Wilson polynomials this is treated in
\cite{29}. In \cite{30} limit cases in the $q$-Askey scheme are also
considered.

Finally we should mention the manuscript \cite{31}.
Here the author extended the duality \cite[(4.2)]{33} for
continuous $q$-ultraspherical polynomials to Macdonald polynomials
and thus obtained the so-called Pieri formula \cite[\S VI.6]{32}
for these polynomials.
\paragraph{Acknowledgement}
I thank Prof.\ M.~A. Pathan and Prof.\ S.~A. Ali for the invitation
to deliver the 2017 R.~P. Agarwal Memorial Lecture and for their
cordiality during my trip to India on this occasion.

\begin{small}
\begin{quote}
T. H. Koornwinder, Korteweg-de Vries Institute, University of
Amsterdam,\\
P.O.\ Box 94248, 1090 GE Amsterdam, The Netherlands;\\
email: {\tt thkmath@xs4all.nl}
\end{quote}
\end{small}

\end{document}